\documentclass{article}
\usepackage{amsmath,amsfonts,amssymb}

  \numberwithin{equation}{section}

\newtheorem{theorem}[equation]{Theorem}
\newtheorem{lemma}[equation]{Lemma}
\newtheorem{corollary}[equation]{Corollary}
\newtheorem{proposition}[equation]{Proposition}
\newtheorem{definition}[equation]{Definition}
\newtheorem{example}[equation]{Example}

\newtheorem{remarks}[equation]{Remarks}

\newcommand{\di}{\operatorname{div}}
\newcommand{\codim}{\operatorname{codim}}

\newcommand{\Pic}{\operatorname{Pic}}
\newcommand{\Cl}{\operatorname{Cl}}
\newcommand{\Ex}{\operatorname{Ext}}
\newcommand{\Exc}{\operatorname{Exc}}
\newcommand{\HOM}{\operatorname{Hom}}
\newcommand{\Sing}{\operatorname{Sing}}
\newcommand{\Supp}{\operatorname{Supp}}
\newcommand{\Jump}{\operatorname{Jump}}
\newcommand{\Div}{{\mathcal D}iv}

\newcommand{\PP}{{\mathbb P}}

\newcommand{\ZZ}{{\mathbb Z}}

\newcommand{\ci}{C^\prime}
\renewcommand{\o}{\omega}
\renewcommand{\O}{{\mathcal O}}
\newcommand{\I}{{\mathcal I}}
\newcommand{\F}{{\mathcal F}}
\newcommand{\E}{{\mathcal E}}
\newcommand{\K}{{\mathcal K}}
\newcommand{\C}{{\mathcal C}}
\newcommand{\G}{{\mathcal G}}
\newcommand{\Hom}{{\mathcal H}om}
\newcommand{\Mod}{{\mathcal M}od}
\newcommand{\Ext}{{\mathcal E}xt}
\newcommand{\Tor}{{\mathcal T}or}

\newenvironment{pf}
{\noindent\textbf{Proof.}} {\hfill{$\square$}\medskip}

{\hfill\textbf{ Q.E.D.}\medskip}

\begin{document}

\title{Linkage on arithmetically Cohen-Macaulay schemes
with application to the classification of curves of maximal
genus.}

\author{Rita Ferraro}

\maketitle
\begin{center}
Dipartimento di Matematica \\
Universit\`a di Roma  Tre\\
Largo San Leonardo Murialdo, 1 - 00146 Roma \\
\end{center}

\section{Introduction}

In algebraic geometry and commutative algebra the notion of
linkage by a complete intersection, which we will here call {\it
classical linkage}, has  been for a long time an interesting and
active topic. In this note we provide a generalization of
classical linkage in  a different context. Namely we will look at
{\it residuals in the scheme theoretic  intersections of
arithmetically Cohen-Macaulay schemes (briefly aCM schemes) of
dimension $r$ (resp. $r+1$) with $r$ hypersurfaces of degree $a_1,
\dots, a_r$}
 (a c.i. of type $(a_1,\dots ,a_r)$
on the aCM scheme, see Def. \ref{c.i}). When the aCM scheme is
singular a c.i. on it may  not be Gorenstein, i.e. its dualizing
sheaf may   not be invertible. If this is the case, classical
linkage, even if suitably generalized, does not apply.

 The main purpose of this article is to prove some results
 related to the invariance of the deficiency module under such
 linkage.
 In the last part of the paper we  show how to apply these
results and techniques to the classification of curves $C$ in
$\PP^n$ of degree $d$ and maximal genus $G(d, n, s)$ among those
not contained in surfaces of degree less than a certain fixed one
$s$. This was the original motivation of this work. A complete
classification theorem has been given for $n=3$ by L. Gruson and
C. Peskine in \cite{gp}, for $n=4$ by L. Chiantini and C.
Ciliberto in \cite{cc} and for $n=5$ by the author in  \cite{f3}.
For $n=3$ and $n=4$ the respective classification Theorems have
been proven with techniques of classical linkage but for $n\geq 5$
this is no longer possible. For $n\geq 5$ and $s\geq 2n-1$ the
classification procedure consists {\it in the precise description
of the linked curve  to $C$ by a certain c.i. on a rational normal
$3$-fold $X$}. In Example \ref{esempio1} we describe this linked
curve in the  easiest  case, i.e. when it is a plane curve. In
Example \ref{esempio2} we construct examples of smooth curves of
maximal genus $G(d, n, s)$ for  every $d$ and $s$ in the range of
Example \ref{esempio1}.

Turning to a detailed presentation of the results, our first one
is the following (see Cor. \ref{cor1}):

\medskip
\noindent{\bf Theorem} \,\,{\it Let $Y_1$, $Y_2$, $Y$ be
projective locally Cohen-Macaulay schemes. If  $Y_1$, $Y_2$ are
geometrically linked by $Y$, or if $Y_1$, $Y_2$ are algebraically
linked by $Y$ and $Y$ is Gorenstein, then}:
\begin{eqnarray*}
\I_{{Y_2}/Y}\otimes\o_Y &\cong&\o_{Y_1}\\
\I_{{Y_1}/Y}\otimes\o_Y &\cong &\o_{Y_2}.
\end{eqnarray*}

 Let now $W\subset\PP^{n-1}$ and $X\subset\PP^n$ be aCM schemes
of dimension $r$ and $r+1$ respectively; throughout the article
$W$ will be often a general hyperplane section of $X$. Let $Z_1$
and $Z_2$ (resp. $Y_1$ and $Y_2$) be the two linked schemes by a
c.i. of type $(a_1, \dots, a_r)$ on $W$ (resp. $X$). Our two main
results are the following isomorphisms of cohomology groups (see
Prop. \ref{linksr}, Th. \ref{linksingr}, Prop. \ref{linksr+1} and
Th. \ref{linksingr+1}) :

\medskip
\noindent{\bf Theorem}
\[
H^0(\I_{{Z_2}/W}\otimes
\o_W(i+ch_W))\cong{H^1(\I_{{Z_1}/W}\otimes\o_W(c-ch_W-i))}^\vee
\]
{\it for $i<\min_j\{a_j\}$, and}
\[
H^1(\I_{{Y_2}/X}\otimes\o_X((i+ch_X))\cong
{H^1(\I_{{Y_1}/X}\otimes\o_X((c-ch_X-i))}^\vee
\]
{\it for every $i$}.

 Here $c=a_1+\cdots+ a_r$, $\o_W$ (resp.
$\o_X$) is the dualizing sheaf of $W$ (resp. of $X$) and  $ch_W$
(resp. $ch_X$) is the smaller integer $k$ such that $\o_W(k)$
(resp. $\o_W(k)$) has sections (see Def. \ref{ch} of {\it
canonical characteristic}).

The first isomorphism above allows us to compute
$h^0(\I_{{Z_2}/W}\otimes \o_W((i+ch_W))$ for low values of $i$ in
terms of the Hilbert function $h_{Z_1}(c-ch_W-i)$ of the residual
scheme $Z_1$. If $Y_1$ is aritmetically Cohen Macaulay, the second
isomorphism implies that $H^1(\I_{{Y_2}/X}\otimes\o_X((i+ch_X))=0$
for every $i$, and therefore the restriction map
$H^0(\I_{{Y_2}/X}\otimes\o_X(i+ch_X)\to
H^0(\I_{{Z_2}/W}\otimes\o_W(i+ch_W)$ is surjective for every $i$
(see Cor. \ref{acm} and Cor. \ref{acmsing}). This means we can
lift curves on $W$ linearly equivalent to $(i+ch_W)H+K_W$ and
passing through a general hyperplane section $Z_2$ of $Y_2$ to
surfaces on $X$ linearly equivalent to $(i+ch_X)H+K_X$ passing
through $Y_2$. Here $H$ denotes the divisor of a hyperplane
section and $K_W$ (resp. $K_X$) is the canonical divisor of $W$
(resp. of $X$).

The technique used to prove the above results allows us to prove
also a formula (see Proposition \ref{genusform}) which relates the
arithmetic genera of the curves $Y_1$ and $Y_2$, linked by a c.i.
$Y$ on the aCM scheme $X$, in the case that $Y$ has no components
contained in the locus where $\o_X$ is not invertible (the {\it
jump locus} of $X$, see Def. \ref{jumplocus}):

\medskip
\noindent{\bf Proposition}
\[
p_a(Y_2)=p_a(Y_1)-p_a(Y)+\deg{K_Y}_{|Y_2}+1.
\]
\medskip

In a first  draft of this  paper the above results were proved in
a -somewhat weaker version- for  linkage by c.i. on rational
normal surfaces and $3$-folds. In the smooth case  the results
were proven using a straightforward generalization of classical
linkage (in particular of Prop. 2.5 of \cite {ps}).
The author warmly thanks the referee who suggests the actual
proofs.  The key point is that  Theorem 1 (see Cor. \ref{cor1})
holds also in the singular case, this permits us to easily
generalize Prop. \ref{linksr} and Prop. \ref{linksr+1} to the
singular case (Th. \ref{linksingr} and Th. \ref{linksingr+1}).

The linkage results presented in this paper, limited to the case
of  c.i. on rational normal $3$-folds in $\PP^5$ and on rational
normal surfaces in $\PP^4$, and the classification for curves of
maximal genus $G(d, n, s)$ in case $n=5$ appeared as part  of my
doctoral dissertation \cite{f1}. The author thanks her advisor
Ciro Ciliberto.

The original version of this paper has been written  while the
author was supported by a INDAM scholarship and the revised
version during a post-doc position of the author at Universit\`a
di Roma Tre.

\section{Preliminaries}

In this section, we collect  the definitions and notation to be
used in this paper, and state some of the basic results of linkage theory.
We introduce the  definitions  of geometric and algebraic
 linkage by a projective scheme $Y$, without
supposing  $Y$ to be a complete intersection.
Moreover we will briefly introduce rational normal scrolls, in particular
what  we need about Weil divisors on them, including linkage.

\begin{definition}
\label{geomlink}
Let $Y_1$, $Y_2$, $Y$ be subschemes of a projective space $\PP$, then
$Y_1$
and $Y_2$ are geometrically linked by $Y$ if:
\begin{enumerate}
\item $Y_1$ and $Y_2$ are equidimensional,  have no embedded components
and  have no common components
\item $Y_1\cup Y_2=Y$, scheme theoretically.
\end{enumerate}
\end{definition}

The following Proposition is essentially Prop. 1.1 of \cite{ps}.
\begin{proposition}
\label{geomisalg}
Let $Y_1$ and $Y_2$ be closed subschemes of $\PP$ geometrically linked by
$Y$, then:
\begin{eqnarray*}
\I_{{Y_1}/Y}&\cong&\Hom_Y (\O_{Y_2}, \O_Y)\\
\I_{{Y_2}/Y}&\cong&\Hom_Y (\O_{Y_1}, \O_Y).
\end{eqnarray*}
\end{proposition}
\begin{pf}
By \cite{ps} Prop. 1.1
 we have that $\I_{{Y_1}/Y}\cong \Hom_\PP (\O_{Y_2}, \O_Y)$
and $\I_{{Y_2}/Y}\cong \Hom_\PP (\O_{Y_1}, \O_Y)$. Since $Y_1$ and
$Y_2$ are both subschemes of $Y\subset \PP$ these isomorphisms can
be rewritten as in the statement.
\end{pf}

\begin{definition}
\label{alglink}
Let $Y_1$, $Y_2$ be  projective schemes, then $Y_1$
and $Y_2$ are algebraically linked by a projective scheme $Y$ containing
them, if:
\begin{enumerate}
\item $Y_1$ and $Y_2$ are equidimensional and have no embedded
components
\item
\begin{eqnarray*}
\I_{{Y_1}/Y}&\cong &\Hom_Y (\O_{Y_2}, \O_Y)\\
\I_{{Y_2}/Y}&\cong &\Hom_Y (\O_{Y_1}, \O_Y).
\end{eqnarray*}
\end{enumerate}
\end{definition}


\begin{remarks}
\label{rem3} If $Y_1$ and $Y_2$ are geometrically linked by $Y$,
then by Prop. \ref{geomisalg} they are also algebraically linked.
Moreover if $Y_1$ and $Y_2$ are algebraically linked by $Y$ and
have no common components, then they are geometrically linked. See
\cite{m} Prop. 5.2.2 (c).
\end{remarks}

\begin{definition}
Let $Y$ be  a projective scheme. The dualizing sheaf of $Y$ is
\[
\o_Y :=\Ext^c_\PP (\O_Y, \o_\PP)
\]
where $Y\hookrightarrow \PP$ is an embedding of $Y$ in some projective
space $\PP$ and
$c=\codim(Y, \PP)$.
\end{definition}
\begin{definition}
\label{jumplocus} Let $Y$ be a projective locally CM scheme and
let $\F$ be a coherent sheaf in $Y$ of dimension $>0$. We define
the jump locus of $\F$ the closed subscheme $\Jump(\F)$ of $Y$
where $\F$ is not locally free.
\end{definition}
\begin{definition}
\label{ch}
 Let $Y$ be a projective scheme, we define the canonical
characteristic of $Y$  the smallest integer $ch_Y$ such that
$h^0(\o_Y(ch_Y))>0$.
\end{definition}

For a proof of the following Theorem the reader may consult
 \cite{e} Th. 21.15 or for more details \cite{f1} Cor. 1.2.3.

\begin{theorem}
\label{ext} Let $Y$ and $X$ be two equidimensional projective
locally Cohen-Macaulay schemes. Suppose $Y\subset X$ and let
$c^\prime$ be $\codim(Y, X)$. Let $\F$ be a sheaf in $\Mod (Y)$.
Then, for every $j\geq 0$:
\[
\Ext^j_Y (\F, \o_Y)\cong \Ext^{{c^\prime}+j}_X(\F, \o_X).
\]
In particular:
\[
\o_Y\cong\Ext^{c^\prime}_X (\O_Y, \o_X).
\]
\end{theorem}

As a corollary of the previous Theorem we show that, supposing $X$
normal, $\o_X$ is the {\it divisorial sheaf} associated to the
canonical divisor $K_X$ of $X$.  For this purpose we  briefly
recall the notion of divisorial sheaves on a normal scheme $X$,
for details and for a more general point of view the reader may
consult the paper of Hartshorne on generalized divisors,
\cite{h1}, \S 2. On a normal scheme, generalized divisors and Weil
divisors are the same (see \cite{h1}, Prop. 2.7).
\begin{definition}
\label{divfascio} Let $X$ be a normal scheme. Let $D$ be a Weil
divisor on $X$. If $K(X)$ denotes the function field of $X$,
 then the sheaf $\O_X(D)$ defined
for  every open set $U\subset X$ as
\[
\Gamma(U, \O_X(D))=\{{f\in K(X)| \di f +D\geq 0 \quad
 \text{on}\quad  U}\}.
\]
is called the divisorial sheaf of $X$.
\end{definition}
It is known (see \cite{h1} Prop. 2.7 and Prop 2.8) that the group
$\Div(X)$ of divisorial sheaves on $X$ is naturally isomorphic to
the group $\Cl(X)$ of Weil divisors modulo linear equivalence.
Moreover there is an equivalence between reflexive sheaves of rank
one and divisorial sheaves.

The following result is known (see e.g. \cite{komo}, Prop. 5.75).
However for sake of self-containedness we give a  direct proof
which  uses the ideas underlying the present work.
\begin{corollary}
\label{I_B|A} If $X$ is an irreducible  normal locally
Cohen-Macaulay projective scheme of positive dimension $r$ , then
the dualizing sheaf $\o_X$ is a twist of the ideal sheaf of a
divisor in X. In particular $\o_X$ is the divisorial sheaf
$\O_X(K_X)$, associated to the canonical divisor $K_X$.
\end{corollary}
\begin{pf}
We can always find a  complete intersection of dimension $r$ and
of certain multi-degree containing $X$.
Let $A$ be a generic such complete intersection and let $\O_A(f)$
be the dualizing sheaf of $A$. Let $B$ be the residual scheme to
$X$ by $A$; by \cite{h1} Prop. 4.1. we see that $X$ and $B$ are
algebraically linked by $A$. Therefore by Theorem \ref{ext} we
have that:
\[\o_X\cong \Hom_A(\O_X, \o_A)\cong \Hom_A(\O_X,
\O_A)(f)\cong\I_{B/A}(f).
\]
Let $X_B$ be the scheme theoretic intersection $X\cap B$.
 Since $\I_{B/A}(f)\cong\o_X$ is supported on $X$,
 tensoring by $\O_X$ do not affect the inclusion $\I_{B/A}(f)\hookrightarrow \O_A(f)$. Therefore
  the exact sequence
$0\to \I_{B/A}(f)\to \O_A(f)\to \O_B(f)$ tensorized by $\O_X$
stays exact, and we find that $\o_X$ is the ideal sheaf of $X_B$
twisted by $\O_X(f)$. Hartshorne's Connectedness Theorem
(\cite{e} Th. 18.12) implies
that $X_B$ is pure of codimension $1$ in $X$, i.e. $X_B$ is a
 Weil divisor of $X$. Therefore $\o_X$ is the divisorial sheaf
$\O_X(fH-X_B)$ (where $H$ is  a hyperplane section). Since
divisorial sheaves and Weil divisors do not depend on closed
subschemes of codimension $\geq 2$, and on the smooth part $X_S$
of $X$ the sheaf ${\o_X}_{|X_S}$ is associated to the canonical
divisor $K_{X_S}$ we have that $\o_X=\O_X(K_X)$.
\end{pf}

\begin{corollary}
\label{cor1} If $Y_1, Y_2, Y$ are  projective locally
Cohen-Macaulay schemes  such that $Y_1, Y_2$ are algebraically
linked by $Y$, then there are the following exact sequences:
\begin{eqnarray*}
0\to\K_1\to\I_{{Y_2}/Y}\otimes\o_Y &\to&\o_{Y_1}\to \C_1\to 0\\
0\to\K_2\to\I_{{Y_1}/Y}\otimes\o_Y&\to &\o_{Y_2}\to \C_2\to 0
\end{eqnarray*}
where $\K_1, \K_2, \C_1$ and $\C_2$ are coherent sheaves on $Y$
with supports contained in $\Jump(\o_Y)\cap Y_1\cap Y_2$. If $Y_1$
and $Y_2$ are geometrically linked by $Y$, or if $Y$ is
Gorenstein, then we have the isomorphisms:
\begin{eqnarray*}
\I_{{Y_2}/Y}\otimes\o_Y &\cong&\o_{Y_1}\\
\I_{{Y_1}/Y}\otimes\o_Y &\cong &\o_{Y_2}.
\end{eqnarray*}
\end{corollary}
\begin{pf}
Let $Y_1$ and $Y_2$ be algebraically linked by $Y$. It is
sufficient to  prove only the first of the two exact sequences.
 Since $\I_{{Y_2}/Y}\cong\Hom_Y(\O_{Y_1}, \O_Y)$ by Definition
\ref{alglink} and $\Hom_Y(\O_{Y_1}, \o_Y)\cong
\Ext^c_\PP(\O_{Y_1}, \o_\PP)\cong\o_{Y_1}$ by Theorem \ref{ext}
(where $Y\hookrightarrow \PP$ is an embedding of $Y$ in some projective
space $\PP$ and
$c=\codim(Y, \PP)$), we have a natural map
\begin{equation}
\label {naturalmap}
 \I_{{Y_2}/Y}\otimes\o_Y\to \o_{Y_1}.
\end{equation}
If $\o_Y$ is locally free we have $\Hom_Y(\O_{Y_1},
\O_Y)\otimes\o_{Y}\cong \Hom_Y(\O_{Y_1}, \o_Y)$, i.e.  map
(\ref{naturalmap}) is an isomorphism (which proves the statement
if $Y$ is Gorenstein); therefore the kernel $\K_1$ and  the
cokernel $\C_1$ of (\ref{naturalmap}) have their supports
contained in $Y_1\cap \Jump(\o_Y)$. The map (\ref{naturalmap}) can
fits into the following commutative diagram with exact rows:
\begin{equation}
\label{diagramma1}
\begin{array}{ccccccccc}
\K &\hookrightarrow &\I_{{Y_2}/Y}\otimes\o_Y &\rightarrow & \o_{Y}
&\rightarrow
& {\o_Y}_{|{Y_2}} & \rightarrow & 0 \\
&& \downarrow && || && \downarrow &&\\
0 &\rightarrow &\o_{Y_1} &\rightarrow & \o_Y & \rightarrow
&\Hom(\I_{{Y_1}/Y},\o_Y) & \rightarrow & 0.
\end{array}
\end{equation}
The top row of the diagram is obtained from  the exact short
sequence $0\to \I_{{Y_2}/Y}\to \O_{Y} \to \O_{Y_2}\to 0$ by
tensoring with $\o_Y$. The kernel $\K$ of the map
$\I_{{Y_2}/Y}\otimes\o_Y\to \O_{Y}\otimes\o_Y$ is $\Tor (\O_{Y_2},
\o_Y)$, which is supported on $\Jump(\o_Y)\cap Y_2$, therefore
$\K_1\cong \K$ has support contained in $\Jump Y\cap Y_1\cap Y_2$.
The bottom  row of the diagram is obtained from  the exact short
sequence $0\to \I_{{Y_1}/Y}\to \O_{Y} \to \O_{Y_1}\to 0$ by
applying   the functor $\Hom (\cdot, \o_Y)$. By Snake's Lemma
$\C_1\cong \ker\{{\o_Y}_{|{Y_2}}\to\Hom(\I_{{Y_1}/Y},\o_Y)\}$,
therefore also $\C_1$ is supported  in $\Jump Y\cap Y_1\cap Y_2$.

Let $Y_1$ and $Y_2$ be geometrically linked by $Y$. It is
sufficient to prove only the first of the two isomorphisms. We use
induction on the dimension $k$ of $Y$. If $\dim Y=0$, then $Y_1$
and $Y_2$ are disjoint, therefore by the first part of this proof
we get $\K_1=\C_1=0$. Let us now suppose that the statement holds
when the dimension  is $k-1$. In particular Cor. \ref{cor1} holds
for the generic hyperplane sections $Z_1, Z_2, Z$ of $Y_1, Y_2, Y$
respectively. Let us consider, for every integer $\alpha$, the
following commutative diagram:
\begin{equation*}
\begin{array}{ccccccc}
\I_{{Y_2}/Y}\otimes \o_Y(\alpha-1) &\rightarrow
&\I_{{Y_2}/Y}\otimes\o_Y(\alpha)
&\rightarrow & \I_{{Z_2}/Z}\otimes \o_Z(\alpha-1) & \rightarrow & 0 \\
 \downarrow && \downarrow && \downarrow &&\\
\o_{Y}(\alpha-1) &\hookrightarrow & \o_Y(\alpha) & \rightarrow &
\o_Z(\alpha-1) & \rightarrow & 0.
\end{array}
\end{equation*}
The cokernels of the three vertical maps are
${\o_{Y}}_{|{Y_2}}(\alpha-1)$, ${\o_{Y}}_{|{Y_2}}(\alpha)$ and $
{\o_{Z}}_{|{Z_2}}(\alpha-1)$ respectively. We want to prove that
the kernels are zero. By induction
$\I_{{Z_2}/Z}\otimes\o_Z(\alpha-1)$ is isomorphic to
$\o_{Z_1}(\alpha-1)$, therefore the third vertical map of the
diagram is injective (the kernel $\K_Z(\alpha-1)$ is isomorphic to
the kernel $\K_1(\alpha-1)$ of map (\ref{naturalmap}), that is
zero by induction). By Snake's Lemma we have that $\K(\alpha-1)$
maps surjectively on $\K(\alpha)$ for every $\alpha$, where $\K$
is the kernel of the map $\I_{{Y_2}/Y}\otimes\o_Y\to \o_Y$. This
implies that $h^0(\K(\alpha-1))\geq h^0(\K(\alpha))$ for every
$\alpha$ and this is possible if and only if $h^0(\K(\alpha))=0$
for every $\alpha$ or if $\K$ is supported over a zero-dimensional
scheme.
The former case clearly yields $\K=0$. We show now that $\K=0$
also if we suppose that $\K$ is supported on a zero-dimensional
scheme. For this purpose it is enough to prove that
$h^0(\I_{{Y_2}/Y}\otimes \o_Y(\alpha) )=0$ for some $\alpha$. Let
$\alpha<ch_Y$, so that $h^0(\o_Y(\alpha))=h^0(\o_Y(\alpha-1))=0$.
Therefore: $
h^0(\I_{{Y_2}/Y}\otimes\o_Y(\alpha-1))=h^0(\K(\alpha-1))=h^0(\K(\alpha))=
h^0(\I_{{Y_2}/Y}\otimes\o_Y(\alpha)). $ This implies
$h^0(\I_{{Y_2}/Y}\otimes\o_Y(\alpha-1))
=h^0(\I_{{Y_2}/Y}\otimes\o_Y(\alpha))=0$ for $\alpha<ch_Y$. Since
$\K\cong \K_1$ we have proven that the natural map
(\ref{naturalmap}) is injective.

We want to prove now that it is also surjective. Let us look at
the exact sequence:
\begin{equation}
\label{sequenza}
 0\rightarrow I_{{Y_2}/Y}\otimes\o_Y\rightarrow
\o_Y\rightarrow {\o_Y}_{|Y_2}\rightarrow 0.
\end{equation}
Following the notation of  Cor. \ref{I_B|A}, the sheaf $\o_Y$ is
$\I_{Y_B/Y}(f)$. Therefore
$\I_{{Y_2/Y}}\otimes\o_Y\hookrightarrow\I_{{Y_B}/Y}(f)\hookrightarrow
\O_Y(f)$ is the ideal sheaf in $Y$ of a scheme
which certainly contains $Y_2$, twisted by $\O_Y(f)$. Hence by the
exact sequence (\ref{sequenza}) we deduce that also the third term
${\o_Y}_{|Y_2}$ is an ideal sheaf in $Y_2$ twisted by
$\O_{Y_2}(f)$, hence torsion free in $Y_2$. Therefore
$\C_1=\ker\{{\o_Y}_{|Y_2}\rightarrow \Hom(\I_{{Y_1}/Y}), \o_Y)\}$,
which is supported on $\Jump(\o_Y)\cap Y_1\cap Y_2$, is zero.
\end{pf}
\medskip

\begin{definition}
\label{c.i} Let $X$ be a projective scheme of dimension $r$; let
$a_i\in \mathbb{N}^+$ and let $1\leq k\leq r$. A complete
intersection (c.i. for short) on $X$ of kind $(a_1, \dots, a_r)$
is an equidimensional projective scheme $Y\subset X$ such that
$\codim(Y, X)=k$, which is scheme theoretic intersection of
Cartier divisors $D_i\in |\O_X(a_i)|$ for $i=1, \dots, k$.
\end{definition}

\medskip

We want now to fix some notation about rational normal scrolls
and point out
what we will need in the next sections.
A rational normal scroll $X\subset \PP$ of dimension $r$ and degree
$f$  is the
image of a projective bundle $\PP(\E)\to \PP^1$ over $\PP^1$
through the morphism $j$ defined by the tautological line bundle
$\O_{\PP(\E)}(1)$, where $\E\cong \O_{\PP^1}(a_1)\oplus \cdots \oplus
\O_{\PP^1}(a_{r})$ with $0\leq a_1\leq \cdots \leq a_{r}$ and
$\sum a_i=f=n-r$. If $a_1=\cdots=a_l=0$, $1\leq l <r$, $X$ is
singular and the  vertex $V$ of $X$ has dimension $l-1$. Let us
denote $\PP(\E)=\tilde X$. The morphism  $j: \tilde{X}\to X$ is a
rational resolution of singularities, i.e. $X$ is {\it normal} and
{\it arithmetically Cohen-Macaulay} and $R^ij_*\O_{\tilde X}=0$
for $j>0$. We will call $j:\tilde X\to X$ the {\it canonical
resolution} of $X$.  It is well known that
$\Pic(\tilde{X})=\mathbb{Z}[\tilde{H}]\oplus\mathbb{Z}[\tilde{R}]$,
where $[\tilde{H}]=[\O_{\tilde{X}}(1)]$ is the hyperplane class
and $[\tilde{R}]=[\pi^*\O_{\PP^1}(1)]$ is the class of the fibre
of the  map $\pi: \tilde X\to \PP^1$. The intersection form on $\tilde X$
is determined by the rule:
\[
{\tilde H}^{r}=f \qquad
{\tilde H}^{r-1}\cdot {\tilde R}=1 \qquad
{\tilde H}^{r-2}\cdot {\tilde R}^2=0.
\]
Let us denote with $X_S$
the smooth part of $X$ and with $\Exc(j)$ the exceptional locus of $j$.
Then
$j:\tilde X\setminus \Exc(j) \to X_S$ is an isomorphism.
Let
$H$ and $R$ be the {\it strict images} of $\tilde H$  and $\tilde
R$
respectively (i.e. the {\it scheme theoretic closure}
$\overline {j({\tilde H}_{|j^{-1}X_S})}$
and
$\overline {j({\tilde R}_{|j^{-1}X_S})}$). Then   we have the following
well known result:
\begin{lemma}
\label{Cl}
Let $X\subset \PP^n$ be a  rational normal scroll of
degree $f$  and let $j:\tilde X\to X$ be its canonical resolution.
Let $\operatorname{Cl}(X)$ be the group of Weil divisors on $X$
modulo linear equivalence. Then:
\begin{enumerate}
\item
If $\codim (V, X)>2$, $\operatorname{Cl}(X)\cong\ZZ [H]\oplus\ZZ
[R]$;
\item
If $\codim(V, X)=2$, $H\sim f R$ and
$\operatorname{Cl}(X)\cong\mathbb{Z}[R]$.
\end{enumerate}
\end{lemma}

We recall here from \cite{f2} the definition of {\it proper} and
{\it (integral)
total transform}
of a Weil divisor in $X$.
In the last section (Example \ref{esempio1}) we will use  proper and
integral total transforms together with \cite{f2} Prop.
4.11 to compute
the multiplicity of the vertex $V$ in the  intersection scheme of two
effective divisors on a
rational normal $3$-fold  $X$ with $\codim(V, X)=2$.

\begin{definition}
\label{proper}
Given a prime divisor $D$ on $X$, the  proper transform $\tilde D$ of
$D$ in $\tilde X$ is the scheme theoretic closure $\overline{j^{-1}(D\cap
X_S)}$. The proper transform of any Weil
divisor in
$X$ is then defined by linearity.
\end{definition}
\begin{definition}
\label{totale} Let $\codim(V, X)=2$ and let $D\sim dR$ be an
effective Weil divisor on $X$, divide $d-1=kf+h$ ($k\geq -1$ and
$0\leq h < f$), we define the {integral total transform} of $D$ as
$\tilde X$ is $D^*\sim (k+1)\tilde H -(f-h-1)\tilde R$.
\end{definition}

Let us define  on $X$ the following coherent sheaves for $a,
b\in \mathbb{Z}$:
\begin{definition}
\[
\O_X(a, b):= j_*\O_{\tilde{X}}(a\tilde{H}+b\tilde{R}).
\]
\end{definition}
We will usually write $\O_X(a)$ instead of $\O_X(a, 0)$. Moreover
for every coherent sheaf $\F$ on $X$ we will write $\F(a, b)$
instead of $\F\otimes\O_X(a, b)$. If the scroll $X$ is smooth,
then the sheaves $\O_X(a, b)$ are the invertible sheaves
associated to the Cartier divisors $\sim aH+bR$ while when $X$ is
singular this is no longer true. In this case we have the
following Proposition which is proved in \cite{f2}
 (Cor. 3.10 and Th. 3.17).

\begin{proposition}
\label{divsheaf}
Let $X\subset \PP^n$ be a singular rational normal scroll of degree $f$,
dimension $r$
and vertex $V$, then:
\begin{enumerate}
\item
If $\codim(V, X)>2$ the sheaf $\O_X(a, b)$ is reflexive for every
$a, b\in\ZZ$ and it is the divisorial sheaf associated to  a Weil divisor
$\sim aH+bR$;
\item
If $\codim(V, X)=2$ the sheaf $\O_X(a, b)$ is reflexive for every
$a, b\in\ZZ$ such that $b<f$; in this case the sheaves
$\O_X(a, b)$ with $a+fb=d$ are all isomorphic to the
the divisorial sheaf associated to a  Weil divisor
$\sim dR$;
\end{enumerate}
\end{proposition}

In the hypotheses  of Prop. \ref{divsheaf}, the dualizing sheaf
$\o_X$ of $X$ is (see \cite{s}):
\begin{equation}
\label{dual}
\o_X=j_*\O_{\tilde X}(K_{\tilde X})=\O_X(-r, f-2).
\end{equation}
By Prop. \ref{divsheaf} we see that $\o_X$ is a divisorial sheaf
(see also Cor. \ref{I_B|A}), therefore the {\it canonical divisor}
of $X$ is $K_X\sim -rH+(f-2)R$. The canonical characteristic
$ch_X$ is then:
\begin{equation}
\label{cancharscroll}
ch_X=r.
\end{equation}


The following result is essentially due to Hartshorne
(Linkage of generalized divisors by a complete intersection:
\cite{h1}, Prop. 4.1). He states it for divisors on a complete
intersection but the same proof goes over as well.

\begin{proposition}{(Linkage of divisors)}
\label{linkdiv} Let $D_1$ be an effective Weil divisor on a normal
scheme $X\subset \PP^n$. Let $F\subset\PP^n$ be a hypersurface
containing $D_1$; let $D$ be the Cartier divisor on $X$ defined by
$F$, then the effective divisor $D_2=D-D_1$ is algebraically
linked to $D_1$ by $D$.
\end{proposition}

\section{Linkage by complete intersection  on aCM schemes}

In this section $W$ is an  aCM scheme of dimension $r$ in
$\PP^{n-1}$ and $X$ is an  aCM scheme of dimension $r+1$ in
$\PP^{n}$.  When $W$ is  smooth we consider algebraic linkage,
where the subschemes need not have distinct components. In the
singular case we consider geometric linkage, where Cor. \ref{cor1}
holds . Let us start with the smooth case and prove the following.
\begin{proposition}
\label{linksr} Let $W\subset\PP^{n-1}$ be a smooth aCM scheme of
dimension $r$. Let $Z_1\subset W$ be  a projective $0$-dimensional
locally Cohen-Macaulay scheme. Let $Z_2\subset W$ be a projective
$0$-dimensional scheme. Let $Z=W\cap F_1\cap\dots\cap F_r$ be a
c.i. on $W$ of type $(a_1,\dots, a_r)$ and let $c=a_1+\cdots a_r$.
Assume that $Z_1$ and $Z_2$ are algebraically linked by $Z$. Let
$ch_W$ be the canonical characteristic  of $W$,
 then for $i<\min_j\{a_j\}$:
\begin{equation*}
H^0(\I_{{Z_2}/W}\otimes \o_W(i+ch_W) )\cong {H
^1(\I_{{Z_1}/W}(c-ch_W-i))}^\vee.
\end{equation*}
This implies, in terms of the Hilbert function of $Z_1$:
\begin{equation*}
h^0(\I_{{Z_2}/W}\otimes \o_W(i+ch_W))=\deg
Z_1-h_{Z_1}(c-ch_W-i).
\end{equation*}
\end{proposition}
\begin{pf}
Just for sake of simplicity we put $i+ch_W=c-v$, where $v$ is an
integer defined by the previous identity. Let us consider the
exact sequence $ 0\to \I_{Z/W}\to \I_{{Z_2}/W}\to \I_{{Z_2}/Z}\to
0 $ tensored by $\o_W(c-v)$:
\[
0\to \I_{Z/W}\otimes \o_W(c-v)\to \I_{{Z_2}/W}\otimes \o_W(c-v)\to
\I_{{Z_2}/Z}\otimes \o_W(c-v)\to 0,
\]
where the right term is  by Cor. \ref{cor1} $\o_{Z_1}(-v)$ (since
$\o_Z={\o_W(c)}_{|Z})$. We want to prove that for
$v>\max_j\{c-a_j-ch_W\}$:
\begin{equation}
\label{primoker} H^0(\I_{{Z_2}/W}\otimes \o_W(c-v))\cong
\ker\{H^0(\o_{Z_1}(-v))\to H^1 (\I_{Z/W}\otimes\o_W(c-v)\}),
\end{equation}
which is equivalent to  prove that
$H^0(\I_{Z/W}\otimes\o_W(c-v))=0$ for $v>\max_j\{c-a_j-ch_W\}$.
For this purpose let us consider the Koszul resolution of
$\I_{Z/W}$ in $\Mod(W)$:
\begin{equation}
\label{z} 0\to E_r\to E_{r-1}\to E_{r-2}\to\cdots \to E_2\to
E_1\to \I_{Z/W}\to 0,
\end{equation}
where the $E_i$'s are finite direct sum $\oplus\O_W(\alpha_i)$
with $\alpha_i\in \ZZ$; in particular $E_1=\oplus_j\O_W(-a_j)$ and
$E_r=\O_W(-c)$. Let $H_1=\ker\{E_1\to\I_{Z/W}\}$ and $H_i$ be
defined recursively as  $H_i=\ker \{E_i\to H_{i-1}\}$, $i=2,\dots,
r-2$. Let us look now at (\ref{z}) tensored
 with the invertible sheaf $\o_W(c-v)$ and let us denote
 the sheaf $\F\otimes\o_W(c-v)$, for every sheaf $\F$ in $W$, with $\tilde\F$.
 Since $W$ is aCM  and
 $E_i=\oplus\O_W(\alpha_i)$, we find that  $H^1(\tilde {H_1})\cong H^{r-2}({\tilde{H}}_{r-2})
 \cong H^{r-1}(\tilde {E_r})=0$. Therefore from the short
exact sequence $0\to H_1\to E_1\to \I_{Z/W}\to 0$ we deduce that
$H^0({\tilde{\I}}_{Z/W})=0$ when $H^0(\tilde{E_1})=\oplus_j
H^0(\o_W(c-a_j-v)))=0$, and this happens for $c-a_j-v<h_W$ for all
$j$ (i.e. $v>\max_j\{c-a_j-ch_W\}$), as we claimed.

At this point we want to prove that
\begin{equation}
\label{secondoker} \Ex^{r-1}_W(\I_{{Z_1}/W}(v),
\o_W)\cong\ker\{H^0(\o_{Z_1}(-v))\to H^r (\o_W(-v))\}.
\end{equation}
Let us consider now exact sequence
\[
0\to\I_{{Z_1}/W}(v)\to\O_W(v)\to\O_{Z_1}(v)\to 0.
\]
 We apply the functor $\HOM(\cdot, \o_W)$:
\[
\dots
\Ex^{r-1}_W(\O_W(v),\o_W)\to\Ex^{r-1}_W(\I_{{Z_1}/W}(v),\o_W)\to\Ex^{r}_W(\O_{Z_1}(v),\o_W)\to\Ex^{r}_W(\O_W(v),\o_W)
\dots
\]
By Serre's duality $\Ex^{r-1}_W(\O_W(v),\o_W)\cong
{H^1(\O_W(v))}^\vee=0$, since $W$ is aCM. Then note that
$\Ex^{r}_W(\O_{Z_1}(v),\o_W)\cong
{H^0(\O_{Z_1}(v))}^\vee=H^0(\o_{Z_1}(-v))$, by Serre's duality on
$W$ for the first isomorphism and on $Z_1$ for the second one. By
\cite{h2}, III, 6.3 $\Ex^{r}_W(\O_W(v),\o_W)\cong
{H^r(\o_W(-v))}$, and we have proven (\ref{secondoker}).

Since $H^r(\o_W(-v))$ "functorially"  contains $H^1
(\I_{Z/W}\otimes\o_W(c-v))$, as one can easily check by looking at
the Koszul complex (\ref{z}) tensored by $\o_W(c-v)$,  we deduce
that  the kernels (\ref{primoker}) and (\ref{secondoker}) are
isomorphic. Therefore, applying Serre's duality to $\Ex^{r-1}_W
(\I_{{Z_1}/W}(v), \o_W)$, it follows that for
$v>\max_j\{c-a_j-ch_W\}$:
\[
H^0(\I_{{Z_2}/W}\otimes \o_W(c-v)
)\cong{H^1(\I_{{Z_1}/W}(v))}^\vee
\]
which proves the first part of the statement
after the substitution $v=c-i-ch_W$.

Moreover, since $W$ is aCM, from the exact sequence
\[
0\to\I_{W/\PP}\to\I_{{Z_1}/\PP}\to\I_{{Z_1}/W}\to 0
\]
we get $H^1(\I_{{Z_1}/W}(k))\cong H^1(\I_{{Z_1}/\PP}(k))$ for
every $k$. Moreover
$h^1(\I_{{Z_1}/\PP}(k))=h^0(\O_{Z_1}(k))-h_{Z_1}(k)=
\deg(Z_1)-h_{Z_1}(k)$. This proves the second part of the
statement.
\end{pf}

We note that Proposition \ref{linksr} can be easily proven using a
straightforward generalization of classical linkage (in particular
of Prop. 2.5 of \cite{ps}). Namely the construction through the
mapping cone of a locally free resolution of $\O_{Z_2}$ in
$\Mod(W)$ from a locally free resolution of $\O_{Z_1}$.
Nevertheless we prefer the proof we have given because it can be
generalized to the singular case:

\begin{theorem}
\label{linksingr}
 Let $W\subset\PP^{n-1}$ be a singular aCM scheme
of dimension $r$. Let $Z_1, Z_2, Z\subset W$ be as in the
hypotheses of Prop. \ref{linksr} and assume moreover that $Z_1$
and $Z_2$ are geometrically linked by $Z$, or that $Z_1$
(equivalently $Z_2$) is contained in the smooth part of $W$. Let
$ch_W$ be the canonical characteristic of $W$,
 then for $i<\min_j\{a_j\}$:
\begin{equation*}
H^0(\I_{{Z_2}/W}\otimes \o_W(i+ch_W) )= {H
^1(\I_{{Z_1}/W}(c-ch_W-i))}^\vee.
\end{equation*}
This implies, in terms of the Hilbert function of $Z_1$:
\begin{equation*}
h^0(\I_{{Z_2}/W}\otimes \o_W(i+ch_W))=\deg
Z_1-h_{Z_1}(c-ch_W-i).
\end{equation*}
\end{theorem}
\begin{pf}
 As in the proof of Prop. \ref{linksr} we want to prove the isomorphism (\ref{primoker}),
 for $v>\max_j\{c-a_j-ch_W\}$. Let us tensor the exact sequence $ 0\to \I_{Z/W}\to
\I_{{Z_2}/W}\to \I_{{Z_2}/Z}\to 0 $ with $\o_W(c-v)$ and  obtain
the exact sequence:
\[
0\to \K\to\I_{Z/W}\otimes \o_W(c-v)\to \I_{{Z_2}/W}\otimes
\o_W(c-v)\to \I_{{Z_2}/Z}\otimes \o_W(c-v)\to 0,
\]
where $\K$, which is a quotient of
$\Tor^1(\I_{{Z_2}/Z},\o_W(c-v))$, has support contained in
$\Jump(\o_W)\cap\Jump(\I_{{Z_2}/Z})\subset\Jump(\o_W)\cap Z_2$.
Let $\mathcal A$  be the kernel of the map $\I_{{Z_2}/W}\otimes
\o_W(c-v)\to \I_{{Z_2}/Z}\otimes \o_W(c-v)$. Since
$H^1(\K)=H^2(K)=0$ we have that $H^1
(\I_{Z/W}\otimes\o_W(c-v))\cong H^1(\mathcal{A})$, hence we
consider $\ker\{H^0(\I_{{Z_2}/Z}\otimes\o_{W}(c-v))\to H^1
(\mathcal A)\cong H^1 (\I_{Z/W}\otimes\o_W(c-v))\}$. By Cor.
\ref{cor1} if $Z_1$ and $Z_2$ are geometrically linked, or if
$Z_1$ or $Z_2$ is contained in the smooth part of $W$, we have
that $\I_{{Z_2}/Z}\otimes \o_W(c-v)\cong \o_{Z_1}(-v)$. Since
$H^0(\I_{Z/W}\otimes\o_W(c-v))=0$ implies $H^0(A)=0$, we prove
(\ref{primoker}) if and only if we prove that
$H^0(\I_{Z/W}\otimes\o_W(c-v))=0$ for $v>\max_j\{c-a_j-ch_W\}$.
For this purpose let us consider the Koszul resolution  (\ref{z})
of $\I_{Z/W}$ in $\Mod(W)$ and note that $\Jump (H_i)$ is
contained in $Z$, for all $i=1, \dots , r-2$. Resolution (\ref{z})
tensored with $\o_W(c-v)$ splits in the following diagrams of
exact sequences, where  we denote ${\tilde{H}}_i=H_i\otimes
\o_W(c-v)$ and ${\tilde{E}}_i=E_i\otimes \o_W(c-v)$:
\begin{equation*}
\begin{array}{ccccccccccccc}
  0 & \to & K_1 & \to & {\tilde{H}}_1 & &\longrightarrow & &{\tilde{E}}_1 & \to &
  \I_{Z/W}\otimes\o_W(c-v) & \to & 0\\
  &  &  &  &  & \searrow &     & \nearrow &  &  &  &  &  \\
  &  &  &  &  &          & C_1 &          &  &  &  &  &  \\
  &  &  &  &  & \nearrow &     & \searrow &  &  &  &  &  \\
  &  &  &  & 0&          &     &          & 0&  &  &  &
\end{array}
\end{equation*}
\begin{equation*}
\begin{array}{ccccccccccccc}
  0 & \to & K_i & \to & {\tilde{H}}_i & &\longrightarrow & &{\tilde{E}}_i & \to &
   {\tilde{H}}_{i-1} & \to & 0\\
  &  &  &  &  & \searrow &     & \nearrow &  &  &  &  &  \\
  &  &  &  &  &          & C_i &          &  &  &  &  &  \\
  &  &  &  &  & \nearrow &     & \searrow &  &  &  &  &  \\
  &  &  &  & 0&          &     &          & 0&  &  &  &
\end{array}
\end{equation*}
for $i=2, \dots , r-2$, and
\begin{equation*}
\begin{array}{ccccccccccccc}
  0 & \to & K_{r-1} & \to & {\tilde{E}}_r & &\longrightarrow & &{\tilde{E}}_{r-1} & \to &
   {\tilde{H}}_{r-2} & \to & 0\\
  &  &  &  &  & \searrow &     & \nearrow &  &  &  &  &  \\
  &  &  &  &  &          & C_{r-1} &          &  &  &  &  &  \\
  &  &  &  &  & \nearrow &     & \searrow &  &  &  &  &  \\
  &  &  &  & 0&          &     &          & 0&  &  &  &
\end{array}
\end{equation*}
The supports of the kernels $K_i$ are contained in the
intersection of $\Jump(\o_W)$ with $\Jump (H_i)$, i.e. $\Supp
K_i\subset \Jump (\o_W)\cap Z$. Therefore looking at the exact
sequences above  we deduce that:
\[
H^1(C_1)\cong H^1({\tilde{H}}_1)\cong H^2(C_2)\dots \cong
H^{r-2}({\tilde H}_{r-2})\cong H^{r-1}(C_{r-1})=0.
\]
Therefore as in the smooth case we have that if
$H^0(\tilde{E_1})=\oplus_j H^0(\o_W(c-a_j-v)))=0$, and this
happens for $c-a_j-v<h_W$ for all $j$ (i.e.
$v>\max_j\{c-a_j-ch_W\}$), then $H^0(\I_{Z/W}\otimes\O_W(c-v))=0$
and we have proven (\ref{primoker}).

As in the proof of Prop. \ref{linksr} one proves the isomorphism:
\[
{H^1(\I_{{Z_1}/W}(v))}^\vee\cong\ker\{H^0(\o_{Z_1}(-v))\to H^r
(\o_W(-v))\},
\]
since it does not depend on the smothness of $W$. Looking again at
the resolution (\ref{z}) tensored by $\o_W(c-v)$, one can easily
check, like in the smooth case, that $H^r(\o_W(-v))$ "functorially"
contains $H^1 (\I_{Z/W}\otimes\o_W(c-v))$
and deduce the first part of the statement.
Like in the proof of Prop. \ref{linksr} one deduces the second
one.
\end{pf}
\medskip

In the next proposition we consider the case of a c.i. of type
$(a_1,\dots ,  a_r)$ on a smooth aCM scheme of dimension $r+1$.
The proof is similar to the proof of Prop. \ref{linksr} and
therefore we give just a brief sketch of it.

\begin{proposition}
\label{linksr+1} Let $X\subset\PP^{n}$ be a smooth aCM scheme of
dimension $r+1$. Let $Y_1, Y_2\subset X$ be  projective
equidimensional  schemes of dimension $1$. Assume that $Y_1$ is
locally Cohen-Macaulay.
 Let $Y\subset X$ be a c.i. of type $(a_1,\dots , a_r)$ on $X$
 and let $c=a_1+\cdots a_r$.
 Assume that $Y_1$ and $Y_2$
are algebraically linked by $Y$. Let $ch_X$ be the canonical
characteristic of $X$, then for every $i$:
\[
H^1(\I_{{Y_2}/X}\otimes \o_X(i+ch_X))\cong
{H^1(\I_{{Y_1}/X}(c-ch_X-i))}^{\vee}.
\]
\end{proposition}
\begin{pf}
We put again, as in the proof of the previous  Proposition,
$i+ch_W=c-v$.
 First we want
to prove that:
\begin{equation}^{}
\label{primokerX}
 H^1(\I_{{Y_2}/X}\otimes \o_X(c-v))\cong
\ker\{H^1(\o_{Y_1}(-v))\to H^2 (\I_{Y/X}\otimes\o_X(c-v))\}.
\end{equation}
Looking at the  exact sequence:
\[
0\to \I_{Y/X}\to \I_{{Y_2}/X}\to \I_{{Y_2}/Y}\to 0
\]
tensored by $\o_X(c-v)$ one see, using Cor. \ref{cor1}, that
(\ref{primokerX}) is equivalent to
$H^1(\I_{Y/X}\otimes\o_X(c-v))=0$ and this is easy to prove using
the  resolution of $\I_{Y/X}$ in $\Mod(X)$. Similarly to the proof
of (\ref{secondoker}) one proves that
\begin{equation}
\label{secondokerX} \Ex^{r}_X(\I_{{Y_1}/X}(v),
\o_X)\cong\ker\{H^1(\o_{Y_1}(-v))\to H^{r+1} (\o_X(-v))\}.
\end{equation}
Since $H^{r+1}(\o_X(-v))$ "functorially" contains $H^2
(\I_{Y/X}\otimes\o_X(c-v))$, one deduces that the kernels
(\ref{primokerX}) and (\ref{secondokerX}) are isomorphic.
Therefore applying Serre's duality to $\Ex^{r}_X (\I_{{Y_1}/X}(v),
\o_X)$ it follows that:
\[
H^1(\I_{{Y_2}/X}\otimes \o_X(c-v)
)\cong{H^1(\I_{{Y_1}/X}(v))}^{\vee}.
\]
and we are done by substituting $v=c-i-ch_X$.
\end{pf}

Let us consider now the singular case and prove the following:
\begin{theorem}
\label{linksingr+1}
 Let $X\subset\PP^{n}$ be a singular aCM scheme
of dimension $r+1$. Let $Y_1, Y_2, Y\subset X$ be as in the
hypotheses of Prop. \ref{linksr+1} and assume moreover that $Y_1$
and $Y_2$ are geometrically linked by $Y$. Let $ch_X$ be the
canonical characteristic of $X$, then for every $i$:
\[
H^1(\I_{{Y_2}/X}\otimes \o_X(i+ch_X))\cong
{H^1(\I_{{Y_1}/X}(c-ch_X-i))}^\vee.
\]
\end{theorem}
\begin{pf}
Following the proof of Prop. \ref{linksr+1} step by step, we want first to
prove the isomorphism (\ref{primokerX}).
 We start from the exact sequence:
\[
0\to K\to\I_{Y/X}\otimes \o_X(c-v)\to \I_{{Y_2}/X}\otimes
\o_X(c-v)\to \I_{{Y_2}/Y}\otimes \o_X(c-v)\to 0,
\]
where $K$ has support contained in $\Jump(\o_X)\cap Y$.
Using  Cor. \ref{cor1} and the same techniques used in the proof of
Th. \ref{linksingr}, one see that (\ref{primokerX}) is equivalent
to $H^1(\I_{Y/X}\otimes\o_X(c-v))=0$.
Using the
resolution of $\I_{Y/X}$ in $\Mod(X)$ tensored with $\o_X(c-v)$
one proves that $H^1(\I_{Y/X}\otimes\o_X(c-v))=0$, as in the
smooth case. The isomorphism (\ref{secondokerX})
does not depend on the smoothness of $X$, therefore it can be
proven exactly as in Prop. \ref{linksr+1}.
Using again the resolution of $\I_{Y/X}$ in $\Mod(X)$ tensored with
$\o_X(c-v)$, one proves also in this case that $H^{r+1}(\o_X(-v))$
"functorially"
contains
$H^2 (\I_{Y/X}\otimes\o_X(c-v))$. Therefore
the kernels (\ref{primokerX}) and (\ref{secondokerX}) are isomorphic. We
conclude as in the proof of Prop. \ref{linksr+1}.
\end{pf}

From Prop. \ref{linksr+1} it follows easily that if we suppose
$Y_1$ arithmetically Cohen-Macaulay, then we can lift "canonical"
divisors $\sim K_W+(i+ch_W)H$ on a general hyperplane section
$W\subset\PP^{n-1}$ of $X$ containing the general hyperplane
section $Z_2$ of $Y_2$ to divisors $\sim K_X+(i+ch_X)H$ on $X$
containing $Y_2$. Namely we have the following Corollary:

\begin{corollary}
\label{acm} In the hypotheses  of Proposition \ref{linksr+1}, if
we suppose $Y_1$ to be arithmetically Cohen-Macaulay, then the map
\[
H^0(\I_{{Y_2}/X}\otimes\o_X(i+ch_X))\to H^0(\I_{{Z_2}/W}\otimes
\o_W (i+ch_W))
\]
is surjective for every $i$.
\end{corollary}
\begin{pf}
Since both $Y_1$ and $X$ are arithmetically Cohen-Macaulay we have
that $h^1(\I_{{Y_1}/X}(k))=h^1(\I_{{Y_1}/\PP}(k))=0$
 for every $k$.
By Prop. \ref{linksr+1}  we  have
$h^1(\I_{{Y_2}/X}\otimes\o_X(i+ch_X))=0$ for every $i$; the
statement follows now from the exact sequence
\[
0\to\I_{{Y_2}/X}\otimes\o_X(i-1+ch_X)\to\I_{{Y_2}/X}\otimes\o_X(i+ch_X)\to
\I_{{Z_2}/W}\otimes\o_X(i+ch_X)\to 0
\]
since $\o_W={\o_X(1)}_{|W}$ and $ch_X=ch_W+1$.
\end{pf}

In the next Proposition we prove that even in the singular case
$H^0(\I_{{Y_2}/X}\otimes \o_X(\alpha))$ represents geometrically
Weil divisors of $X$ linearly equivalent to $K_X+\alpha H$
containing $Y_2$. The proof is similar to the one we have used in
Cor. \ref{cor1} to prove the exact sequence (\ref{sequenza}),
therefore we omit it.
\begin{proposition}
In the hypotheses of Theorem \ref{linksingr+1} we have the
following exact sequence:
\[
0\to\I_{{Y_2}/X}\otimes \o_X\to \o_X\to {\o_X}_{|Y_2}\to 0.
\]
\end{proposition}
With this in mind we can state  Cor. \ref{acm}  also in the
singular case:
\begin{corollary}
\label{acmsing} In the hypotheses  of Theorem \ref{linksingr+1},
if we suppose $Y_1$ to be arithmetically Cohen-Macaulay, then the
map
\[
H^0(\I_{{Y_2}/X}\otimes\o_X(i+ch_X))\to H^0(\I_{{Z_2}/W}\otimes
\o_W (i+ch_W))
\]
is surjective for every $i$.
\end{corollary}
\begin{pf}
Following the proof of Cor. \ref{acm} one have to consider the
following exact sequence:
\[
0\to K \to
\I_{{Y_2}/X}\otimes\o_X(-1)\to\I_{{Y_2}/X}\otimes\o_X\to
\I_{{Z_2}/W}\otimes\o_X\to 0
\]
where  kernel $K$
has support contained in $\Jump(\o_X)\cap Z_2$. Let $H$ be the
kernel of the map $\I_{{Y_2}/X}\otimes\o_X\to
\I_{{Z_2}/W}\otimes\o_X$. Since
$h^1(H(i+ch_X))=h^1(\I_{{Y_2}/X}\otimes\o_X(i-1+ch_X))=0$ for all
$i$, we are done.
\end{pf}

 The next result is a formula which
relates the arithmetic genera of the  curves $Y_1$, $Y_2$ and $Y$.

\begin{proposition}
\label{genus} In the hypoteses of Prop. \ref{linksr+1} or of
Theorem \ref{linksingr+1}, if we suppose that  $Y$ has no
 components contained in $\Jump(\o_X)$, we have the following
formula, relating the arithmetic genera of the linked curves:
\begin{equation}
\label{genusform} p_a(Y_2)=p_a(Y_1)-p_a(Y)+\deg({K_Y}_{|Y_2})+1.
\end{equation}
\end{proposition}
\begin{pf}
First we tensor by
 $\o_Y$  the exact sequence
\[
0\to\I_{{Y_2}/Y}\to \O_Y\to\O_{Y_2}\to 0.
\]
 By Cor. \ref{cor1}  this
gives the exact sequence (see exact sequence (\ref{sequenza})):
\begin{equation}
\label{seq1} 0\to\o_{Y_1}\to \o_Y\to{\o_Y}_{|Y_2}\to 0.
\end{equation}
With the notation of Cor. \ref{I_B|A} we know that $\o_Y$ is
$\I_{{Y_B}/Y}(f+c)$, where $Y_B$ is, in this case,  the  intersection
of $X_B$ with the $r$ hypersurfaces which cut $Y$ in $X$. Let us
remark that $X_B$ passes through the jump locus  of $X$, otherwise
$X$ would be there locally complete intersection, hence
Gorenstein,  and this can not happen by definition of jump locus.
If some components  of $Y$ are contained in $\Jump(\o_X)$, then $Y_B$
contains such components and it is clearly not a divisor of $Y$. On the
other
hand,  since we can always arrange
the complete intersection $A$ of Cor. \ref{I_B|A} in such a way that $X_B$
does not contain any components of $Y$ outside $\Jump(\o_X)$,
if $Y$ has no components contained in $\Jump(\o_X)$
the scheme $Y_B=X_B\cap Y$ is a divisor in $Y$.
In this case  $\o_Y$ is the divisorial sheaf in $Y$ associated to $K_Y$
and  ${\o_Y}_{|Y_2}$ is a divisorial sheaf in $Y_2$,
associated to the divisor $D={K_Y}_{|Y_2}$, i.e.
$\O_{Y_2}(D)={\o_Y}_{|Y_2}$. By Riemann-Roch theorem (for
eventually singular curves) we then obtain
\begin{equation}
\label{rel1}
\chi(\O_{Y_2}(D))=\chi({\o_Y}_{|Y_2})=1-p_a(Y_2)+\deg({K_Y}_{|Y_2}).
\end{equation}
 Formula (\ref{genusform}) follows now by (\ref{rel1}) and
(\ref{seq1}).
\end{pf}

\section{An application to the classification of curves of
maximal genus}

In this section we will show  some examples of application
of the techniques developed in the previous sections to the classification
of curves of maximal genus $G(d, n,  s)$ in $\PP^n$.
Let us first  summarize some results of
\cite{ccd} and some other preliminary  facts
useful to introduce
the problem.

From now on, let $C$ be an integral, nondegenerate curve of degree
$d$ and arithmetic genus $p_a(C)$ in $\PP^n$, with
$d>\frac{2s}{n-2}\Pi_{i=1}^{n-2} {((n-1)!)}^ {\frac{1}{n-1-i}}$
and $s\geq n-1$ (later we will assume $s\geq 2n-1$). Assume $C$ is
not contained on surfaces of degree $<s$ and define $m, \epsilon,
w, v, k, \delta$ as follows:

divide $d-1=sm+\epsilon$, $0\le \epsilon\le s-1$ and
$s-1=(n-2)w+v$, $v=0, \dots , n-3$;

if $\epsilon< w(n-1-v)$, divide $\epsilon=kw+\delta$, $0\le\delta<w$;

if $\epsilon \ge w(n-1-v)$, divide $\epsilon+n-2-v=k(w+1)+\delta$, $0\le
\delta<w+1$.

It is a result of \cite{ccd} (section 5) that the genus $p_a(C)$ is
bounded by the function:
\[
G(d, n, s)=1+\frac{d}{2}(m+w-2)-\frac{m+1}{2}(w-3)+\frac{vm}{2}(w+1)+\rho
\]
where $\rho=\frac{-\delta}{2}(w-\delta)$ if $\epsilon<w(n-1-v)$ and
$\rho=\frac{\epsilon}{2}-\frac{w}{2}(n-2-v)-\frac{\delta}{2}(w-\delta+1)$
if $\epsilon\ge w(n-1-v)$.

If $Z$ is a general hyperplane section of $C$ and $h_Z$ is the Hilbert
function of $Z$, then the difference $\Delta h_Z$ must be bigger than the
function $\Delta h$ defined by:
\[
\Delta h(r)=
\begin{cases}
0 & \text{if} \quad r<0 \cr (n-2)r+1 & \text{if} \quad  0\le r \le
w \cr s & \text{if} \quad  w< r \le m \cr s+k-(n-2)(r-m) &
\text{if} \quad m< r \le m+\delta \cr s+k-(n-2)(r-m)-1 & \text{if}
\quad m+\delta < r \le m+w+e \cr 0 & \text{if} \quad r>m+w+e
\end{cases}
\]
where $e=0$ if $\epsilon<w(n-1-v)$ and $e=1$ otherwise (\cite{ccd} Prop.
0.1).
\begin{proposition}
\label{extcurve}
If $p_a(C)=G(d, n, s)$, then $C$ is arithmetically Cohen-Macaulay and
 $\Delta h_Z(r)=\Delta h(r)$ for all $r$.
Moreover $Z$ is contained on a reduced curve $\Gamma$ of degree
$s$ and maximal genus $G(s, n-1)={w\choose 2}+wv$ in $\PP^{n-1}$
(Castelnuovo curve). $\Gamma$ is unique and, when we move the
hyperplane, all these curves $\Gamma$'s patch togheter giving a
surface $S\subset \PP^n$ of degree $s$ through $C$ (Castelnuovo
surface).
\end{proposition}
\begin{pf}
See \cite{ccd} Prop. 6.1, Prop. 6.2 and Cor. 6.3.
\end{pf}

We recall that a {\it Castelnuovo curve} in $\PP^n$ is a
nondegenerate reduced and irreducible curve of degree $d$ and
maximal arithmetic genus $G(d, n)$. Castelnuovo in 1893 found the
bound $G(d, n)$ and he went on to give a complete geometric
description of those curves which achieved his bound, they lie on
surfaces of minimal degree. The reader may consult \cite{ha1} for
all the details. By {\it Castelnuovo surface} we mean  a
nondegenerate reduced and irreducible surface in $\PP^n$ whose
general hyperplane section is a Castelnuovo curve in $\PP^{n-1}$
(see \cite{ha2}).

\begin{proposition}
\label{surf}
The surface $S$ of Prop. \ref{extcurve}
is irreducible and if $s\geq 2n-1$ it  lies on a
 rational normal  $3$-fold $X\subset\PP^n$. As a divisor
on $X$ the surface $S$ is linearly equivalent to $(w+1)H-(n-3-v)R$
(or $wH+R$ if $v=0$). If $n=6$ and if $s$ is even  there is the
further possibility that the surface $S$ lies in a cone over the
Veronese surface in $\PP^5$ and is the complete intersection with
a hypersuface not containing the vertex.
 \end{proposition}
\begin{pf}
$S$ is irreducible since $C$ is irreducible and is not contained on
surfaces of degree $<s$.
The rest of the statement follows using the characterization of
Castelnuovo surfaces given in \cite{ha2}.

\end{pf}

\begin{proposition}
\label{m+1}
There exists a hypersurface
$F_{m+1}$ of degree $m+1$, passing through $C$ and not containing $S$.
\end{proposition}
\begin{pf}
For a general hyperplane section $\Gamma$ of $S$, the Hilbert
function $h_\Gamma$ is known (see e.g. \cite{ha1} Th. 3.7); in
particular we have $\Delta h_\Gamma (r)=\Delta h_Z (r)$ when $0\le
r\le m$ and hence $h^0({\I}_{C/\PP} (r))=h^0({\I}_{S/\PP}(r))$
when $0\le r\le m$. For $r=m+1$ one shows $\Delta h_\Gamma
(m+1)<\Delta h_Z (m+1)$ and this implies $h^0({\I}_{C/{\PP}}
(m+1))>h^0({\I}_{S/{\PP}}(m+1))$.
\end{pf}

\medskip
Let us suppose $s\geq 2n-1$. By the  Prop. \ref{surf} a  curve
$C\subset\PP^n$ of maximal genus $G(d, n, s)$ lies then on a
rational normal $3$-fold $X$ (except in the case where $S$ lies in
a cone over a Veronese surface, which we do  not intend to go
through). Let $F_{w+1}$ be a hypersurface of degree $w+1$ cutting
out on $X$ the surface $S$. By Prop. \ref{m+1} we can consider  on
$S$ the curve $C^\prime$
 residual to $C$ by the intersection
with the hypersurface $F_{m+1}$. Since $\deg(C^\prime)<\deg(C)$,
then $C^\prime$ does not contain $C$. Choosing in $X$ a
sufficiently general  divisor  $D\sim (n-3-v)R$ (or $D\sim H-R$ in
case $S\sim wH+R$) linked to $S$ by  $X\cap F_{w+1}$, then the
residual  scheme on $X$ to $C$  by the c.i. $X\cap F_{w+1}\cap
F_{m+1}$ is a curve which we call $C^{\prime\prime}$. When $v=n-3$
then $S=X\cap F_{w+1}$ and of course $C^\prime=C^{\prime\prime}$;
otherwise $C^{\prime\prime}$ is the union of $C^\prime$ with a
curve $C_D$ contained in $D$, therefore $C_D$ is formed by $n-3-v$
distinct plane curves of degree $m+1$ or, in case $S\sim wH+R$,
$C_D$ is the complete intersection on  $D\sim H-R$ by a
hypersurface of degree $m+1$. Letting $Z^\prime$,
$Z^{\prime\prime}\subset W$ be general hyperplane sections of
$C^\prime$ and $C^{\prime\prime}$ respectively, we have the
following Lemma:

\begin{lemma}
\label{first} Let $C^{\prime\prime}\subset X$ be as in the
previous notation. Then for $i\leq w, m$
\[
h^0(\I_{C^{\prime\prime}/X}(i, n-4))\geq
h^0(\I_{Z^{\prime\prime}/W}(i, n-4)) = \sum_{r=m+w-i+1}^\infty
\Delta h(r).
\]
Moreover if $h^0(\I_{Z^{\prime\prime}/W}(i-1, n-4))=0$ and
$h^0(\I_{Z^{\prime\prime}/W}(i, n-4))=h>0$, then
$h^0(\I_{C^{\prime\prime}/X}(i-1, n-4))=0$ and
$h^0(\I_{C^{\prime\prime}/X}(i, n-4))=h$.
\end{lemma}
\begin{pf}
$C$ and $C^{\prime\prime}$ are geometrically linked by $Y=X\cap
F_{w+1}\cap F_{m+1}$ since they are equidimensional, have no
common components ($C$ is irreducible and $C^\prime$ does not
contain $C$) and  no embedded components ($Y$ is arithmetically
Cohen Macaulay). By Prop.  \ref{linksr} (if $X$ is smooth) or Th.
\ref{linksingr+1} (if $X$ is singular) we know that $
h^0(\I_{Z^{\prime\prime}/W}(i, n-4))= d-h_Z(m+w-i) $. Then note
that for every $k$ we have $d-h_Z(k)=d+\Delta h_Z(k+1)-
h_Z(k+1)=d+\sum_{r=k+1}^t \Delta h_Z(r)-h_Z(t)
=\sum_{r=k+1}^\infty \Delta h_Z(r)$ because for $t$ big we have
$h_Z(t)=d$, and that, by Prop. \ref{extcurve}, $\Delta
h_Z(r)=\Delta h (r)$ for all $r$.  Since $\o_X\cong\O_X( -3, n-4)$
and $\o_W\cong\O_W(-2, n-4)$, by Cor. \ref{acm} we have that
\begin{equation}
\label{bigger} h^0(\I_{C^{\prime\prime}/X}(i, n-4))\geq
h^0(\I_{Z^{\prime\prime}/W}(i, n-4)).
\end{equation}
Let us consider the exact sequence
\[
0\to K(k+3)\to\I_{C^{\prime\prime}/X}(k-1, n-4)\to
\I_{C^{\prime\prime}/X}(k, n-4)\to \I_{Z^{\prime\prime}/W}(k,
n-4)\to 0,
\] where $K$ is the kernel of the map $\I_{C^{\prime\prime}/X}(-4,
n-4)\to \I_{C^{\prime\prime}/X}(-3, n-4)$ and is supported on
$\Sing(X)\cap Z^{\prime\prime}$. If, for $k=i-1$, we have
$h^0(\I_{Z^{\prime\prime}/W}(i-1, n-4))=0$, since
$h^1(K(i+2))=h^0$ we have that  $h^0(\I_{C^{\prime\prime}/X}(i-2,
n-4))\geq h^0(\I_{C^{\prime\prime}/X}(i-1, n-4))$, which implies
$h^0(\I_{C^{\prime\prime}/X}(i-1, n-4))=0$. In this hypothesis,
for $k=i$, we have an injection $H^0(\I_{C^{\prime\prime}/X}(i,
n-4)) \hookrightarrow H^0(\I_{Z^{\prime\prime}/W}(i, n-4))$ and
therefore by (\ref{bigger})
 $h^0(\I_{C^{\prime\prime}/X}(i,
n-4)) = h^0(\I_{Z^{\prime\prime}/W}(i, n-4))$.
\end{pf}

We rewrite now the genus formula (\ref{genusform}) for linked
curves contained in a rational normal three-fold  $X$:
\begin{proposition}
\label{genere} Let $Y_1$ and $Y_2$ be two
$1$-dimensional projective schemes contained in a rational normal
$3$-fold $X\subset \PP^n$ which is smooth or whose vertex is a
point. Let $Y_1$ be locally Cohen-Macaulay. Assume that $Y_1$ and
$Y_2$ are geometrically linked  by a complete intersection
$Y=X\cap F_a\cap F_b$ of type $(a, b)$ on $X$ (if $X$ is smooth it
is enough to suppose $Y_1$ and $Y_2$ algebraically linked by $Y$).
Then:
\begin{equation}
\label{generescroll}
 p_a(Y_2)=p_a(Y_1)-p_a(Y)+(a+b-3)\cdot
\deg(Y_2)+(n-4)\cdot\deg(R_{|{Y_2}})+1.
\end{equation}
\end{proposition}
\begin{pf}
Since $\Jump(\o_X)$ is a point, we can apply Proposition
\ref{genus}. Formula (\ref{generescroll}) follows from formula
(\ref{genusform}) since
${K_Y}_{|Y_2}\sim{(a+b-3)H}_{|Y_2}+{(n-4)R}_{|Y_2}$.
\end{pf}

The strategy is to classify all the curves of maximal genus in $\PP^n$ for
arbitrary $n$ by classifing
the  linked curves $C^\prime$'s. A complete classification
Theorem when $n=4$ is proved in \cite{cc} and  when $n=5$
in \cite{f1} (and in the forthcoming work \cite{f3}). Depending on the numerical
parameters ( $\epsilon$, $w$, $v$, $k$) associated to $C$  and
on the
type of the scroll $X$ the analysis goes on case by case.
In the following example we want to show the simplest non trivial case in
the classification
procedure, when $C^\prime$ is a plane curve (the trivial case is
$C^\prime=\emptyset$).   It should be remarked that while in $\PP^3$ the curve
$C^\prime$ is always degenerate this is no longer true for $n\geq 4$
(see \cite{cc} and \cite{f1} for $n=4, 5$).

\begin{example}
\label{esempio1}
Let $s\geq 2n-1$ and let $d>\frac{2s}{n-2}\Pi_{i=1}^{n-2} {((n-1)!)}^
{\frac{1}{n-1-i}}$; divide $s-1=(n-2)w+v$, $v=0, \dots
, n-3$ and divide $d-1=sm+\epsilon$, $0\leq\epsilon\leq s-1$.
Suppose $s-2-w\leq\epsilon\leq s-2$. Let $C\subset\PP^n$ be a curve
of maximal genus $G(d, n, s)$. Then  the linked curve $C^\prime$
 is a plane curve
 of degree $s-\epsilon-1$.
In case that the vertex of $X$ is  a line, $C^\prime$ will not contain
this line
as a component.
\end{example}

Here  we suppose for the sake of  simplicity  that $v=n-3$ (the
result can be proved with similar arguments
for every $v$), i.e. we
put
ourselves in the simplest case $C^\prime=C^{\prime\prime}$; with this
assumption we always have $e=1$, i.e.
$\epsilon\geq w(n-1-v)=(n-3)(w+1)$, hence we write
  $\epsilon+1=k(w+1)+\delta$ with
$k=n-3$ and $\delta\leq w$. In this case $C$ and $C^\prime$ are
(geometrically) linked  by a c.i. $Y=X\cap F_{w+1}\cap F_{m+1}$ on
$X$.  Applying Lemma \ref{first} for $i=0$ (of course
$h^0(\I_{C^{\prime}/X}(-1, n-4))=0$) we compute:
\[
h^0(\I_{C^{\prime}/X}(0, n-4))=
n-4.
\]
Let us exclude for now the case $X$ singular along a line. The
linear system $|\O_X(0, n-4)|$ is composed with a rational pencil,
i.e. we have $\pi: X\to \PP^1$ and $|\I_{C^\prime/X}(0,
n-4))|=\pi^* \G$, where $\G$ is a linear subsystem of
$|\O_{\PP^1}(n-4)|$. Since $h^0(\O_{\PP^1}(n-4))=n-3$, this
implies that $|\I_{C^{\prime}/X}(0, n-4)|$ has a fixed part; in
this case, since $h^0(\O_X(0, a))=a+1$ for every $a\geq 0$, the
fixed part of $|\I_{C^\prime}(0, n-4)|$ is $\sim R$ and the moving
part is equal to the whole $|\O_X(0, n-5)|$. Therefore we conclude
that $C^\prime$ is contained in a plane $\pi\sim R$.


Let us consider now the case when the vertex of $X$ is a line. We
want to conclude as in the previous case that
$|\I_{{C^\prime/X}}(0, n-4)|$ has a fixed part. So let us suppose
that $|\I_{C^{\prime}/X}(0, n-4)|$ has no fixed part, which
implies that the support of $C^\prime$ is the singular line of
$X$. By Bertini's Theorem  the generic divisor in the
corresponding linear subsystem $\G$ of $\PP^1$ is union of $n-4$
distinct points in a rational normal  curve $C_{n-2}$ of degree
$n-2$ (a $(n-2)$-plane section of $X$), which span a $\PP^{n-5}$.
Therefore we can choose a basis $\{D_1, \dots , D_{n-4}\}$ in the
linear system $|\I_{C^{\prime}/X}(0, n-4)|$ such that $D_i$ is
union of $n-4$ distinct planes of $X$ for every $i$ and such that
the linear space spanned by  each $D_i$ is $<D_i>\cong\PP^{n-3}$
and $D_i=X\cap <D_i>$. In this situation the base locus of
$|\I_{C^{\prime}/X}(0, n-4)|$, which is equal to $D_1\cap \cdots
\cap D_{n-4}=X\cap <D_1>\cap \cdots \cap <D_{n-4}>$, is
necessarily the singular line $l\cong \PP^1$ of $X$ counted with
multiplicity one, but this is not possible since $|\I_{l|X}(0,
n-4))|= |\O_X(0, n-4)|$ and we  have a contradiction. Therefore,
as in the previous case, we conclude that $C^\prime$ is contained
in a plane $\pi\sim R$. We claim now that $C^{\prime}$ cannot
contain the singular line of $X$ as a component. In fact in this
case both $S$ and $F_{m+1}$ would pass through it and their proper
transforms $\tilde S$ and $\tilde {F}_{m+1}$ on the canonical
resolution $\tilde X$ of $X$ would be $\tilde S\sim (w+1-a)\tilde
H+ (n-2)a\tilde R$ and $\tilde {F}_{m+1}\sim (m+1-b)\tilde
H+(n-2)b\tilde R$ with $a, b\geq 1$. In this case, since $C$ is
irreducible (therefore it does not contain the singular line),
$C^\prime$ would contain the singular line with multiplicity
$\alpha$ which we compute using \cite{f2} Prop. 4.11 as:
$\alpha=S^*\cdot F^*_{m+1}\cdot\tilde H-\tilde S\cdot \tilde
{F}_{m+1}\cdot \tilde H =(w+1)\tilde H\cdot(m+1){\tilde
H}^2-\tilde S\cdot \tilde {F}_{m+1}\cdot \tilde H= ab(n-2)$, where
$S^*$ and $F^*_{m+1}$ are respectively the integral total
transform of $S$ and $F_{m+1}$ (Def. \ref{totale}). But since
$C^\prime$ is contained in a plane $\pi\sim R$ by the same kind of
computation we conclude that $C^\prime$ would contain the singular
line with multiplicity $\beta=S^*\cdot R^*\cdot\tilde H-\tilde
S\cdot \tilde R\cdot \tilde H = (w+1)\tilde H\cdot (\tilde H
-(n-3)\tilde R)\cdot\tilde H-((w+1-a)\tilde H+(n-2)a\tilde R)\cdot
\tilde R\cdot \tilde H=a$ and this is in contradiction with the
previous value.

\bigskip
In the next example we show that in the  case of Example
\ref{esempio1} {\it smooth} curves of maximal genus do always
exist. Moreover we explicitly construct such curves  on a smooth
rational normal $3$-fold.  It is interesting to note that it is
not always possible to construct curves of maximal genus on a {\it
smooth} rational normal $3$-fold. There are cases (for some values
of $d$ and $s$) where the construction is possible only on  a
rational normal $3$-fold whose vertex is a point and where  genus
formula (\ref{generescroll}) holds, as showed in
 \cite{f3} (Prop. 4.2 case 6) and Theorem 5.2 case $k=v=1$)
 for $n=5$.
 The existence of
curves of maximal genus in $\PP^5$ is proved for all cases in \cite{f3}.
We state first the following, easy to prove,
result (see \cite{ro} Lemma 1 pg. 133) which we will use later.

\begin{lemma}
\label{rogora}
Let $X$ be a smooth $3$-fold. Let $\Sigma$ be a linear system of surfaces
of $X$ and let $\gamma$ be a curve contained in the base locus of $\Sigma$.
Suppose that the generic surface of $\Sigma$ is smooth at the generic
point of $\gamma$ and that it has at least a singular point which is variable
in $\gamma$. Then all the surfaces of $\Sigma$ are tangent along $\gamma$.
\end{lemma}

\begin{example}
\label{esempio2}
For every $d$ and $s$ in the range
 of $Example$ \ref{esempio1} there exists a smooth curve
$C\subset \PP^n$ of maximal genus $G(d, n, s)$.
\end{example}
 For the sake of simplicity we treat only the cases $v=n-3, n-4$, i.e.
$s=(n-2)(w+1)$
and
$s=(n-2)w+n-3$. The other cases can be treated in a similar way.

Let us suppose  $v=n-3$.
Let $X\subset \PP^n$ be a smooth rational normal $3$-fold of degree
$n-2$ and
let $\pi\sim R$ be a plane contained in $X$.
Let $D$ be a smooth curve on $\pi$ of degree \,
$0\leq \deg D=w+1-s+\epsilon+1=\epsilon+1-(n-3)(w+1)\leq w$
(possibly $D=\emptyset$). If we consider the union of $D$
with any plane curve $C^\prime\subset\pi$ of degree
$w+1-\deg D=s-\epsilon-1$, then
 there exists a hypersurface $F_{w+1}$
of degree $w+1$ cutting out on $\pi$ the union $C^\prime\cup D$.
Therefore the linear system $|\I_{{D}/X}(w+1)|$ of divisors on $X$
cut by hypersurfaces of degree $w+1$ through $D$ is not empty and
cut on $\pi$ the linear system $D+|\O_\pi(s-\epsilon-1)|$.
Moreover the linear system $|\I_{{D}/X}(w+1)|$ contains
 the linear
subsystem $L+|\O_X(w)|$, where $L$ is a fixed hyperplane section
containing $\pi$, that has fixed part $L$ and no other base
points. This implies that $D$ is  the base locus of all
$|\I_{{D}/X}(w+1)|$ and that $|\I_{{D}/X}(w+1)|$ is not composed
with a pencil, because in this case every element in the system
would be a sum of algebraically equivalent divisors, while the
divisors in  $L+|\O_X(w)|$ are obviously not of this type.
 By Bertini's
Theorem we can then  conclude that the generic divisor in
$|\I_{{D}/X}(w+1)|$ is an irreducible surface $S$ of degree
$(n-2)(w+1)=s$   smooth outside $D$. We claim that  $S$ is in fact
smooth at every point $p$ of $D$. To see this, by Lemma
\ref{rogora}, it is enough to prove that, for every $p\in D$,
there exists a surface in $|\I_{D/X}(w+1)|$ which is smooth at
$p$, and that for a generic point $q\in D$, there exist two
surfaces in $|\I_{D/X}(w+1)|$ with distinct tangent planes at $q$.
In fact, for every $p\in D$ we can always find a surface $T$ in
the linear system $|\O_X(w)|$ which does not pass through $p$,
therefore the surface $L+T$ is smooth at $p$ with tangent plane
$\pi$. Moreover a generic surface in the linear system
$|\I_{D/X}|$ which cut $D$ on $\pi$ has at $p$ tangent plane
$T_p\neq \pi$.

Let $\ci\subset\pi$ be the linked curve to $D$ by the intersection
$\pi\cap S$. Let us consider the linear system
$|\I_{{C^\prime}/S}(m+1)|$ of divisors cut on $S$ by the
hypersurfaces of degree $m+1$ passing through $C^\prime$. With the
same argument used above we conclude that
 this linear system is not composed with a pencil, it has $C^\prime$
as a  fixed part
 and no other base points.
Therefore by Bertini's theorem we deduce that the generic  curve
$C=S\cap F_{m+1}-C^\prime$ in the movable part of
the linear system is irreducible, smooth
and has
the required degree
$d=s(m+1)-s+\epsilon+1$.
By Clebsch formula one computes:
\[
p_a(C^\prime)=\frac{1}{2}((n-2)w+n-4-\epsilon)((n-2)w+n-5-\epsilon).
\]
Moreover $\deg(R\cap C^\prime)=0$. Substituting these expressions
in the genus formula (\ref{generescroll}) we find that $p_a(C)$
has the maximal value $G(d, n, s)$, therefore $C$ is the required
curve.

We consider now the case $v=n-4$.
 Let $\pi\sim R$  and  $p\sim R$ be two distinct planes contained in $X$.
Let $D$  be a smooth curve on $\pi$ of degree \, $0\leq \deg D=
\epsilon+2-(n-3)(w+1)\leq w$ (possibly $D=\emptyset$). Let us
consider the linear system $|\I_{{D\cup p}/X}(w+1)|$
 of divisors on $X$ cut by
hypersurfaces of degree $w+1$ containing the plane $p$ and passing
through $D$. This linear system is not empty since hypersurfaces
which are union of a hyperplane containing the plane $p$ and of a
hypersurface of degree $w$ passing through $D$ cut on $X$ divisors
in the system. From this description one can see that $|\I_{{D\cup
p}/X}(w+1)|$ is not composed with a pencil and that its
 base locus is  $p\cup D$. By Bertini's Theorem the generic element
in the movable part of the linear system is an irreducible surface
$S\sim (w+1)H-R$ of degree $s$, smooth outside $D$. By the same
argument used in the previous case we can prove that $S\cup p$ is
smooth at every point of $D$, but since $D\cap p =\emptyset$ this
means that $S$ is smooth at $D$. Let $\ci\subset\pi$ be the linked
curve to $D$ by the intersection $S\cap \pi$. Let us consider the
linear system $|\I_{{C^\prime}/S}(m+1)|$, which is not empty since
$\deg\ci <m+1$ and has base locus equal to the curve $\ci$.
As in the previous case we deduce that the generic  curve
$C=S\cap F_{m+1}-C^\prime$ in the movable part of
this linear system is irreducible, smooth
and has
the required degree $d=s(m+1)-\deg(C^\prime)$.
By generality the hypersurface $F_{m+1}$ does not contain the plane $p$
and cut on it a curve $C_1$ of degree $m+1$.
Let
$C^{\prime\prime}=C^\prime\cup C_1$; by construction the curve
$C^{\prime\prime}$
is
 linked to $C$ by a c.i. on $X$ of type $(w+1, m+1)$

By Noether's formula one computes:
\[
p_a(C^{\prime\prime})=\frac{1}{2}((n-2)w+n-5-\epsilon)((n-2)w+n-6-\epsilon)
+\frac{1}{2} m(m-1)-1.
\]
Moreover $\deg(R\cap C^{\prime\prime})=0$. Substituting these
expressions in the genus formula (\ref{generescroll}) we find that
$p_a(C)$ has the maximal value $G(d, n, s)$. Therefore $C$ is the
required curve.

\end{document}